\documentclass[10pt,journal,onecolumn]{IEEEtran}

\usepackage{graphicx, bm, amsmath, siunitx, longtable,tabularx, float, mathtools, booktabs, eqnarray, color, amsfonts, amssymb, graphics, setspace, makecell}

\usepackage{cite}  

\usepackage[version=4]{mhchem}
\pdfoutput=1

\newcommand{\B}[1]{{\bm #1}}
\newcommand{\T}{^{\mbox{\tiny T}}}

\newcommand{\dd}{\; \text{d}}

\newcommand{\ces}{constrained expressions}
\renewcommand{\ce}{constrained expression}

\usepackage{etoolbox}

\begin{document}

\title{Least-squares Solutions of Eighth-order Boundary Value Problems using the Theory of Functional Connections$^*$\thanks{$^*$This work was supported by a NASA Space Technology Research Fellowship, Johnston [NSTRF 2019] Grant \#: 80NSSC19K1149 and Leake [NSTRF 2019] Grant \#: 80NSSC19K1152}}

\author{Hunter Johnston\thanks{Aerospace Engineering, Texas A\&M University, College Station, TX, 77843 USA (hunterjohnston@tamu.edu, leakec@tamu.edu, mortari@tamu.edu)}, \and Carl Leake, and \and Daniele Mortari}

\maketitle

\begin{abstract}
    This paper shows how to obtain highly accurate solutions of eighth-order boundary-value problems of linear and nonlinear ordinary differential equations. The presented method is based on the Theory of Functional Connections, and is solved in two steps. First, the Theory of Functional Connections analytically embeds the differential equation constraints into a candidate function (called a {\it constrained expression}) that contains a function that the user is free to choose. This expression always satisfies the constraints, no matter what the free function is. Second, the free-function is expanded as a linear combination of orthogonal basis functions with unknown coefficients. The constrained expression (and its derivatives) are then substituted into the eighth-order differential equation, transforming the problem into an unconstrained optimization problem where the coefficients in the linear combination of orthogonal basis functions are the optimization parameters. These parameters are then found by linear/nonlinear least-squares. The solution obtained from this method is a highly accurate analytical approximation of the true solution. Comparisons with alternative methods appearing in literature validate the proposed approach.
\end{abstract}

\section{Introduction}

This paper has been motivated by several articles dedicated to the solution of high-order boundary-values problems (BVPs) including, fourth-order \cite{R01}, sixth-order \cite{R02}, eighth-order \cite{R00,R03,R04,R05,R06,R07,R08,R09,R10,R11,R12,R13,R14,R15,R16,R17,R18}, $2m$-order \cite{R19}, and higher-order \cite{R20,R21,R22,R23} BVPs. This paper focuses specifically on the eighth-order BVPs because of the volume of research done on them, which is covered in Refs. \cite{R00,R03,R04,R05,R06,R07,R08,R09,R10,R11,R12,R13,R14,R15,R16,R17,R18}. These references list many existing scientific problems requiring solutions of high-degree BVPs. For example, eighth-order BVPs appear in the physics of specific hydrodynamic stability problems (infinite horizontal layer of fluid heated from below and under rotation) when instability sets in as overstability \cite{stability_book}, and in orthotropic cylindrical shells under line load \cite{cylindraical_shells}. From a theoretical point of view, the study of the existence and uniqueness of high order boundary value problems is presented in Ref. \cite{math_BVP_book} and studied further in Ref. \cite{existence_8th_order}. 

The technique presented in this paper is rooted in functional interpolation expressions. These expressions are particularly well suited to solve differential equations. This has been shown in Ref. \cite{U-TFC}, the seminal paper on the Theory of Functional Connections (TFC), and in subsequent articles, showing its application to solve linear \cite{LDE} and nonlinear ODEs \cite{NDE}. The TFC formalized the method of analytical constraint embedding (a.k.a. functional interpolation), since it provides expressions representing all functions satisfying a set of specified constraints. 

The general equation to derive these interpolating expressions, named \ces, follows as,
\begin{equation}\label{eq:gen_ce}
    y(x,g(x)) = g(x) + \sum_{k=1}^n \eta_k s_k(x)
\end{equation}
where $g (x)$ is the free function, $\eta_k$ are unknown coefficients to be solved by imposing the $n$ constraint conditions, and $s_k (x)$ are ``support functions,'' which are a set of $n$ linearly independent functions. In prior work \cite{LDE,NDE} as well as in this paper, the $s_k (x)$ support function set has been selected as the monomial set. 

The $\eta_k$ coefficients are computed by imposing the constraints using Eq. (\ref{eq:gen_ce}). Then, once the expressions of the $\eta_k$ coefficients are obtained, they are back substituted into Eq. (\ref{eq:gen_ce}) to produce the {\it constrained expression}, a functional representing all possible functions satisfying the specified set of constraints. The use of this \ce\ has already been applied to many areas of study, including solving low-order differential equations \cite{LDE, NDE}, hybrid systems \cite{hybrid_tfc}, and optimal control problems, including energy-optimal landing, energy-optimal intercept \cite{EOL_EOI}, and fuel-optimal landing \cite{FOL}. Furthermore, this technique has been successfully used to embed constraints into machine learning frameworks \cite{SVM,DeepTfc}, in quadratic and nonlinear programming \cite{QP_NLP}, and in a variety of other applications \cite{Selected}. In addition, this technique has been generalized to $n$-dimensions \cite{M-TFC,M-TFC-PDE}, providing functionals representing all possible $n$-dimensional manifolds subject to constraints on the value and arbitrary order derivative of $n-1$ dimensional manifolds.

\section{Derivation of the Constrained Expression for Eighth-order Boundary-value Problems}\label{sec:SolveBVP}
In this paper, we consider the solution of an eighth-order BVP via the TFC. In general, the problem can be posed in its implicit form as,
\begin{equation}\label{eq:gen_8th_de}
    F(x,y,y',\hdots,y^{(8)}) = 0 \quad \text{subject to:} \begin{cases} y^{(k)}(x_0) = y^{(k)}_i \\  
    y^{(k)}(x_f) = y^{(k)}_f \\ \end{cases} \quad \text{for} \quad k = 0,1,2,3
\end{equation}
where the notation $y^{(k)}(x) := \frac{\dd^{k} y(x)}{\dd x^{k}}$ is used to denote the $k^{th}$ derivative of $y(x)$ with respect to $x$. Now, in order to embed the eight constraints, we can set $n = 8$ in Eq. (\ref{eq:gen_ce}) leading to the expression,
\begin{equation}\label{eq:candidate}
    y(x,g(x)) = g(x) + \B{\eta}\T \B{s}(x)
\end{equation}
where
\begin{equation*}
    \B{\eta} = \begin{Bmatrix} \eta_1, & \eta_2, & \eta_3, & \eta_4, & \eta_5, & \eta_6, & \eta_7, & \eta_8 \end{Bmatrix}\T
\end{equation*}
\begin{equation*}
    \B{s}(x) = \begin{Bmatrix} 1, & x, & x^2, & x^3, & x^4, & x^5, & x^6, & x^7 \end{Bmatrix}\T
\end{equation*}
Now, according to the theory developed in Ref. \cite{U-TFC}, a system of equations can be constructed by evaluating the candidate function defined by Eq. (\ref{eq:candidate}) at the constraint locations and setting the function equal to the specified constraint value. For example, the constraint on the function at the initial value (i.e. $y(x_i) = y_i$) is applied as such,
\begin{equation*}
    y_i = y(x_i,g(x_i)) = g_i + \eta_1 + \eta_2 x_i + \eta_3 x_i^2 + \eta_4 x_i^3 + \eta_5 x_i^4 + \eta_6 x_i^5 + \eta_7 x_i^6 + \eta_8 x_i^7.
\end{equation*}
This can be done for the remaining seven constraint conditions, and the resulting system of equations can be expressed in a compact form,
\begin{equation*}
    \begin{Bmatrix} 
    y_i - g_i \\ y_f - g_f \\
    y^{(1)}_i - g^{(1)}_i \\ y^{(1)}_f - g^{(1)}_f \\
    y^{(2)}_i - g^{(2)}_i \\ y^{(2)}_f - g^{(2)}_f \\
    y^{(3)}_i - g^{(3)}_i \\ y^{(3)}_f - g^{(3)}_f
    \end{Bmatrix} = \begin{bmatrix} 
    1 & x_i & x_i^2 & x_i^3 & x_i^4 & x_i^5 & x_i^6 & x_i^7 \\
    1 & x_f & x_f^2 & x_f^3 & x_f^4 & x_f^5 & x_f^6 & x_f^7 \\
    
    0 & 1 & 2 x_i & 3x_i^2 & 4x_i^3 & 5x_i^4 & 6x_i^5 & 7x_i^6 \\
    0 & 1 & 2 x_f & 3x_f^2 & 4x_f^3 & 5x_f^4 & 6x_f^5 & 7x_f^6 \\
    
    0 & 0 & 2 & 6x_i & 12x_i^2 & 20x_i^3 & 30x_i^4 & 42x_i^5 \\
    0 & 0 & 2 & 6x_f & 12x_f^2 & 20x_f^3 & 30x_f^4 & 42x_f^5 \\
    
    0 & 0 & 0 & 6 & 24x_i & 60x_i^2 & 120x_i^3 & 210x_i^4 \\
    0 & 0 & 0 & 6 & 24x_f & 60x_f^2 & 120x_f^3 & 210x_f^4 \\
    \end{bmatrix} \B{\eta}.
\end{equation*}
This system of equations can be solved for the unknown $\B{\eta}$ coefficients and organized in the form,
\begin{equation}\label{eq:ce_solved}
\begin{aligned}
    y(x,g(x)) = g(x) &+ 
    \beta_1(x) (y_i - g_i) + \beta_2(x) (y_f - g_f) \\ &+ 
    \beta_3(x) \left( y^{(1)}_i - g^{(1)}_i\right) + 
    \beta_4(x) \left( y^{(1)}_f - g^{(1)}_f\right)  +
    \beta_5(x) \left( y^{(2)}_i - g^{(2)}_i \right) \\ &+ 
    \beta_6(x) \left( y^{(2)}_f - g^{(2)}_f \right) + 
    \beta_7(x) \left( y^{(3)}_i - g^{(3)}_i \right) + 
    \beta_8(x)\left( y^{(3)}_f - g^{(3)}_f \right),
\end{aligned}
\end{equation}
where the $\beta_k(x)$ terms, called switching functions, are solely a function of the independent variable. This technique is general for any domain $x \in [x_0, x_f]$, for example the general expression for $\beta_1(x)$ is,
\begin{equation*}
\begin{aligned}
    \beta_1(x) = \frac{1}{(x_f-x_i)^7}\Big[(x-x_f)^4 \Big(20 x^3 &-7 x_i \left(10 x^2+4 x x_f+x_f^2\right)+10 x^2 x_f \\ &+4 x x_f^2+21 x_i^2 (4 x+x_f)+x_f^3-35 x_i^3\Big)\Big]
\end{aligned}
\end{equation*}
However, for ease of presentation, since all problems presented in this paper, except for Problem \#5, are defined on the domain $x \in [0,1]$, we will express these switching functions in terms of this integration range. For completeness, the support functions for two general points $x_i$ and $x_f$ (i.e. the switching functions for Eq. \eqref{eq:ce_solved}) are provided in appendix \ref{app:SupportFunctions}. The $\beta_k$ terms for $x \in [0,1]$ are summarized below in Eqs. (\ref{eq:beta1}-\ref{eq:beta8}).
\begin{eqnarray}
\beta_1(x) &=& 20 x^7-70 x^6+84 x^5-35 x^4+1\label{eq:beta1}\\
\beta_2(x) &=& -20 x^7+70 x^6-84 x^5+35 x^4 \label{eq:beta2}\\
\beta_3(x) &=& 10 x^7-36 x^6+45 x^5-20 x^4+x \label{eq:beta3}\\
\beta_4(x) &=& 10 x^7-34 x^6+39 x^5-15 x^4 \label{eq:beta4}\\
\beta_5(x) &=& 2 x^7-\dfrac{15 x^6}{2}+10 x^5-5 x^4+\dfrac{x^2}{2} \label{eq:beta5}\\
\beta_6(x) &=& -2 x^7+\dfrac{13 x^6}{2}-7 x^5+\dfrac{5 x^4}{2} \label{eq:beta6}\\
\beta_7(x) &=& \dfrac{x^7}{6}-\dfrac{2 x^6}{3}+x^5-\dfrac{2 x^4}{3}+\dfrac{x^3}{6} \label{eq:beta7}\\
\beta_8(x) &=& \dfrac{x^7}{6}-\dfrac{x^6}{2}+\dfrac{x^5}{2}-\dfrac{x^4}{6}\label{eq:beta8}
\end{eqnarray}
With the solution of the $\beta_k(x)$ terms, the \ce\ is fully solved and represents all possible functions satisfying the boundary-value constraints. More specifically, by substituting the \ce\ and its derivatives into the differential equation it is transformed into a differential equation of $g(x)$ and its derivatives subject to no constraints. Therefore, the differential equation expressed in Eq. (\ref{eq:gen_8th_de}) can now be rewritten such that,
\begin{equation}\label{eq:DE_tilde_g}
    \tilde{F}(x,g,g',\hdots,g^{(8)}) = 0
\end{equation}
To solve this differential equation, prior work \cite{LDE,NDE,EOL_EOI,FOL,Selected} has expanded $g(x)$ in terms of some known basis, such that,
\begin{equation*}
    g(x) = \B{\xi}\T \B{h}(z)
\end{equation*}
where $z = z(x)$, $\B{\xi}$ is an $m \times 1$ vector of unknown coefficients, and $\B{h}(z)$ is an $m \times 1$ vector of some known basis (in this paper Chebyshev polynomials are used). Attention must be paid to the terms used in the basis expansion. The basis functions in $g(x)$ \emph{must} be linearly independent of the support functions $s_k(x)$ in order to solve the system via least-squares. If any of the terms in $g(x)$ is not linearly independent of the support functions, then the matrix that must be inverted when performing least squares will be singular. Thus, the terms that are not linearly independent of the support functions must be skipped in the expansion of $\B{h}(z)$. In this problem, our support functions span from the monomial term $x^0$ to $x^7$; therefore, the Chebyshev polynomial expansion must start at the eighth-order term. Furthermore, in general, the basis functions may not be defined on the same range as the problem domain (i.e. Chebyshev and Legendre polynomials are defined on $z \in [-1, +1]$, Fourier series is defined on $z \in [-\pi, +\pi]$, etc.). Therefore, the basis domain ($z$) must be mapped to the problem domain ($x$), which can be done using the simple linear equations,
\begin{equation*}
    z = z_0 + \frac{z_f-z_0}{x_f-x_0}(x - x_0) \quad \longleftrightarrow \quad x = x_0 + \frac{x_f-x_0}{z_f-z_0}(z - z_0).
\end{equation*}
Furthermore, all subsequent derivatives of the free-function $g(x)$ are defined as,
\begin{equation*}
    \frac{\dd^{n} g}{\dd x^{n}} = \B{\xi} \T  \frac{\dd^{n} \B{h}(z)}{\dd z^{n}} \left(\frac{\dd z}{\dd x}\right)^{n},
\end{equation*}
where by defining,
\begin{equation}\label{eq:c_mapping}
c := \frac{\dd z}{\dd x} = \frac{z_f - z_0}{x_f - x_0}
\end{equation}
the expression can be simplified to, 
\begin{equation}\label{eq:basis_derv_n}
    \frac{\dd^{n} g}{\dd x^{n}} = c^{n} \B{\xi} \T  \frac{\dd^{n} \B{h}(z)}{\dd z^{n}} = c^{n} \B{\xi} \T \B{h}^{(n)}(z).
\end{equation}
This defines all mappings from the basis domain into the problem domain. With the expression of $g(x)$ in terms of a known basis, we can rewrite the \ce\ given in Eq. (\ref{eq:ce_solved}) in the form,
\begin{equation}\label{eq:general_ce_a_b}
    y(x,\B{\xi}) = a(x,\B{\xi}) +  b(x)
\end{equation}
where $a(x,\B{\xi})$ is a function that is zero where the constraints are defined and $b(x)$ is a function which equals the constraint where the constraints are defined. In fact, if $g(x)$ was selected such that $g(x)=0$ (meaning $\B{\xi} = \B{0}$) the expression would simplify to an interpolating function $y(x,\B{0}) =  b(x)$. Furthermore, since the $\B{\xi}$ vector shows up linearly in the expression of $a(x,\B{\xi})$, Eq. (\ref{eq:general_ce_a_b}) can also be written as,
\begin{equation}\label{eq:general_ce_a_b_xi}
    y(x,\B{\xi}) = \B{a}(x)\T \B{\xi} +  b(x),
\end{equation}
where $a(x,\B{\xi})$ now becomes a vector equation $\B{a}(x)$. This can be seen by expanding Eq. (\ref{eq:general_ce_a_b_xi}),
\small
\begin{align*}
    &y(x,\B{\xi}) = \\
    &\overbrace{\Big(\B{h}(z) 
    - \beta_1(x) \B{h}_i
    - \beta_2(x) \B{h}_f
    - \beta_3(x) \B{h}^{(1)}_i 
    - \beta_4(x) \B{h}^{(1)}_f 
    - \beta_5(x) \B{h}^{(2)}_i
    - \beta_6(x) \B{h}^{(2)}_f
    - \beta_7(x) \B{h}^{(3)}_i
    - \beta_8(x) \B{h}^{(3)}_f \Big)}^{\B{a(x)}} \phantom{}\T \B{\xi} \\
   &+ \underbrace{\beta_1(x) y_i
    + \beta_2(x) y_f
    + \beta_3(x) y^{(1)}_i 
    + \beta_4(x) y^{(1)}_f 
    + \beta_5(x) y^{(2)}_i
    + \beta_6(x) y^{(2)}_f
    + \beta_7(x) y^{(3)}_i
    + \beta_8(x) y^{(3)}_f}_{b(x)}.
\end{align*}
\normalsize
The subsequent derivatives follow by simply taking the derivatives of the $\B{h}(z)$ and $\beta_k(x)$ terms. That is, the form of subsequent derivatives of the \ce\ remains the same and we can generally write the \ce\ up to the eighth-order derivative as shown in Eq. (\ref{eq:general_ce_a_b}),
\begin{equation}\label{eq:general_ce_a_b_derv}
    \begin{cases}
    y^{(1)}(x,\B{\xi}) = \B{a}^{(1)}(x)\T \B{\xi} + b^{(1)}(x)\\
    y^{(2)}(x,\B{\xi}) = \B{a}^{(2)}(x)\T \B{\xi} + b^{(2)}(x)\\
    y^{(3)}(x,\B{\xi}) = \B{a}^{(3)}(x)\T \B{\xi} + b^{(3)}(x)\\
    y^{(4)}(x,\B{\xi}) = \B{a}^{(4)}(x)\T \B{\xi} + b^{(4)}(x)\\
    y^{(5)}(x,\B{\xi}) = \B{a}^{(5)}(x)\T \B{\xi} + b^{(5)}(x)\\
    y^{(6)}(x,\B{\xi}) = \B{a}^{(6)}(x)\T \B{\xi} + b^{(6)}(x)\\
    y^{(7)}(x,\B{\xi}) = \B{a}^{(7)}(x)\T \B{\xi} + b^{(7)}(x)\\
    y^{(8)}(x,\B{\xi}) = \B{a}^{(8)}(x)\T \B{\xi} + b^{(8)}(x)\\
    \end{cases}
\end{equation}
where $\B{a}^{(n)}(x) := \frac{\dd \B{a}(x)}{\dd x^n}$ is the derivative of the $\B{a}(x)$ function; the $\B{a}(x)$ function also includes the derivative of $\B{h}(z)$, which follows Eq. (\ref{eq:basis_derv_n}) such that $\frac{\dd \B{h}^{(n)}(z)}{\dd x^n} = c^{n}\B{h}^{(n)}(z)$ where $c$ is defined in Eq. (\ref{eq:c_mapping}). With this adjustment to the \ce, the transformed differential equation defined in Eq. \ref{eq:DE_tilde_g} can be reduces to a function of only $x$ and the unknown vector $\B{\xi}$,
\begin{equation}\label{eq:simple_de}
    \tilde{F}(x,\B{\xi}) = 0,
\end{equation}
which may be linear or nonlinear in the unknown parameter $\B{\xi}$. Lastly, in order to solve this equation numerically, we must discretize the domain into $N+1$ points. Since in this paper we consider the linear basis $\B{h}(z)$ as Chebyshev orthogonal polynomials, the optimal distribution of $N+1$ points is provided by Chebyshev-Gauss-Lobatto
collocation points \cite{Colloc,ChebCol}, defined as,
\begin{equation*}
z_k = -\cos\left(\frac{k \pi}{N}\right) \quad \text{for} \quad k=0,1,\hdots, N,
\end{equation*}
and the map from $z \rightarrow x$ has been previously defined. As compared to the uniform point distribution, the collocation point distribution allows a much slower increase of the condition number as the number of basis functions, $m$, increases. In general, we can define the residual of our differential equation in Eq. (\ref{eq:simple_de}) for each discretized point,
\begin{equation}\label{eq:de_discrete}
\tilde{F}(x_k,\B{\xi}) = 0.
\end{equation}
For a linear differential equation $F$ (and therefore a linear differential equation $\tilde{F}$) the \ce\ and its derivatives will show up linearly and therefore will remain linear in the unknown $\B{\xi}$ term. This leads to the form,
\begin{equation*}
    \mathbb{A}(\B{x})\B{\xi} + \B{b}(\B{x}) = \B{0}
\end{equation*}
where the matrix $\mathbb{A}$ is composed of a linear combination of $\B{a}(x)$ and its derivatives discretized over $x_k$ where $\B{x} = [x_0,\cdots,x_k,\cdots,x_N]\T$. Furthermore, it follows that $\B{b}(\B{x})$ composed of a linear combination of $b(x)$, its derivatives, and a potential forcing term $f(x)$ for the discrete values of $x$. This system can now be easily solved with any available least-squares technique. All numerical solutions in this paper utilize a scaled QR method to perform the least-squares.

In the case of a nonlinear differential equation, Eq. (\ref{eq:de_discrete}) can be expressed as a loss function at each discretization point,
\begin{equation*}
    \mathcal{L}(\B{\xi}_i) = \begin{bmatrix}  \tilde{F}(x_0,\B{\xi}) \\ \vdots \\ \tilde{F}(x_k,\B{\xi}) \\ \vdots \\  \tilde{F}(x_N,\B{\xi}) \end{bmatrix}_{\B{\xi} = \B{\xi}_i}
\end{equation*}
and the system can be solved by an iterative least-squares method where the Jacobian is defined as,
\begin{equation*}
    \mathcal{J}(\B{\xi}_i) = \begin{bmatrix} \dfrac{\partial \tilde{F}(x_0,\B{\xi})}{\partial \B{\xi}} \\ \vdots \\ \dfrac{\partial \tilde{F}(x_k,\B{\xi})}{\partial \B{\xi}} \\ \vdots \\ \dfrac{\partial \tilde{F}(x_N,\B{\xi})}{\partial \B{\xi}} \end{bmatrix}_{\B{\xi} = \B{\xi}_i}
\end{equation*}
where $\B{\xi}_i$ represents the current step's estimated $\B{\xi}$ parameter. The parameter update is provided by,
\begin{equation*}
    \B{\xi}_{i+1} = \B{\xi}_i - \Delta \B{\xi}_i
\end{equation*}
where the $\Delta \B{\xi}_i$ can be defined using classic least-squares,
\begin{equation*}
    \Delta \B{\xi}_i = \Big(\mathcal{J}(\B{\xi}_i)\T \mathcal{J}(\B{\xi}_i) \Big)^{-1} \mathcal{J}(\B{\xi}_i)\T \mathcal{L}(\B{\xi}_i)
\end{equation*}
or in this paper through a QR decomposition method. This process is repeated until either the absolute value of the loss function is below some tolerance $\epsilon$, or until or the the $L_2$ norm of the loss function continues to increases, which is specified by the following conditions,
\begin{equation*}
    L_2 [\mathcal{L}(\B{\xi}_i)] < \varepsilon \qquad \text{or} \qquad
    L_2 [\mathcal{L}(\B{\xi}_{i+1})] > L_2 [\mathcal{L}(\B{\xi}_{i})].
\end{equation*}
In this paper, this tolerance was set as twice the value of machine-level precision for double point precision, $\epsilon = 4.4409 \times 10^{-16}$

\section{Parameter Initialization for Nonlinear Problems}\label{sec:init}
\begin{figure}[H]
    \centering\includegraphics[width=.65\linewidth]{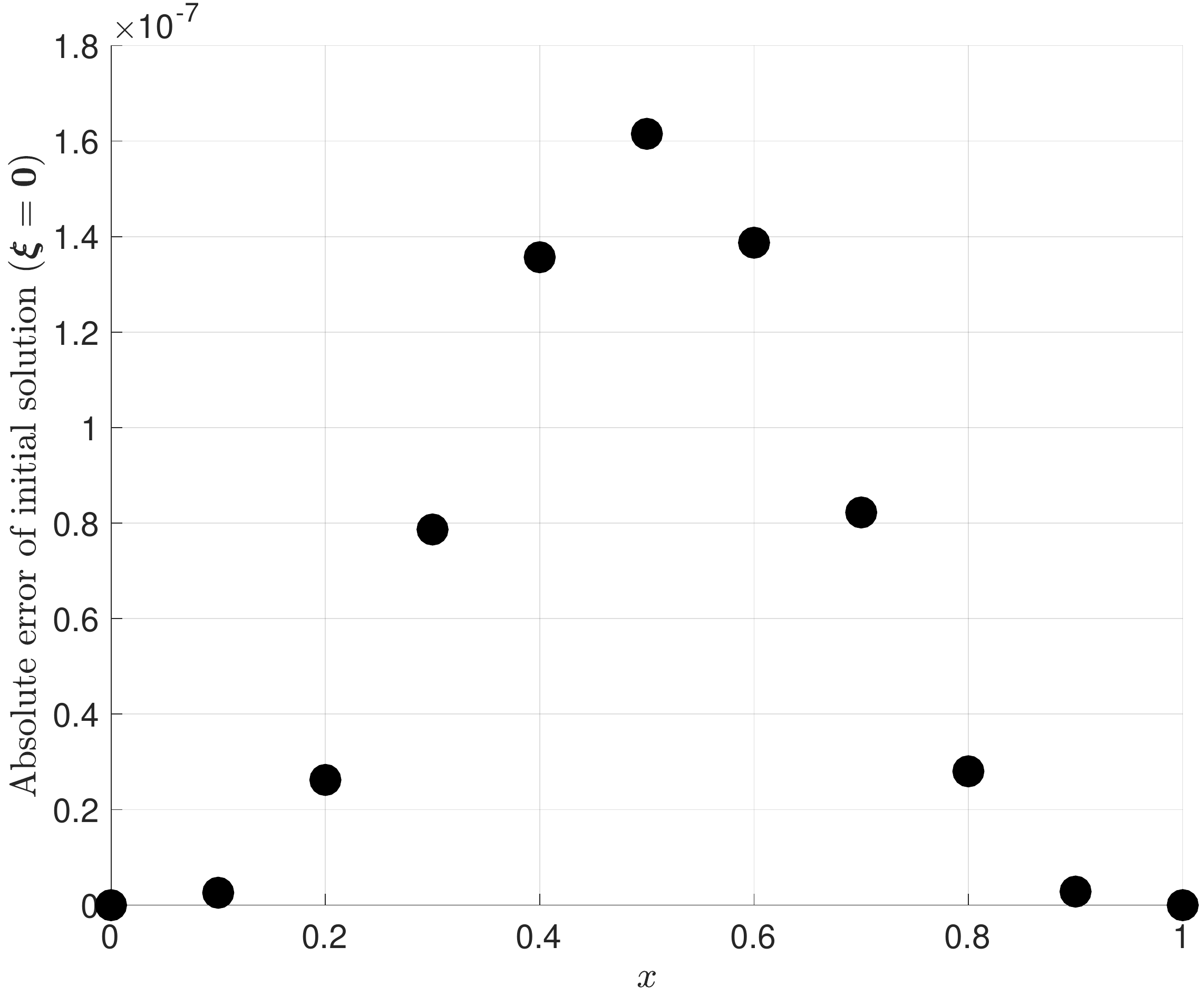}
    \caption{Initialization error of the solution of Problem \#4 by imposing $\B{\xi}_0 = \B{0}$.}
    \label{fig:init}
\end{figure}

For nonlinear problems the $\B{\xi}$ must be initialized at the beginning of the iterative least-squares process. In this paper, the initialization was chosen to be $\B{\xi}_0 = \B{0}$ for all nonlinear problems. Setting the coefficient vector equal to zero is synonymous with selecting $g(x) = 0$, or in other words, choosing the \ce\ with the simplest interpolating polynomial satisfying all of the problem constraints. In the case of BVPs, the solution lies somewhere around this initial guess. Although introduced and solved in a later section, consider Problem \#4 which involves solving the differential equation,
\begin{equation*}
    y^{(8)}(x) + y^{(3)}(x)\sin(y(x)) = e^{x}(1+\sin(e^x)) \quad x \in[0,1].
\end{equation*}
This problem is highlighted in this section because it had the largest initialization error of the three nonlinear differential equations presented. Figure \ref{fig:init} displays the error of the solution due to the initialization technique of  $\B{\xi} = \B{0}$. It can be seen in this figure, the error is on the order of $10^{-7}$ for this specific case. The iterative least-squares will then reduce this error to close to machine-level precision. An astute reader will notice that the TFC method at initialization produces a more accurate solution than the techniques developed in \cite{R00,R17}. A more detailed explanation of this result is discussed in the conclusion.

\section{Numerical Solution}
This section compares the TFC method with competing methods on a variety of problems, both linear and nonlinear. For each problem, the differential equation and boundary conditions are presented, followed by a table that compares the absolute error of the two methods on a grid of 11 equidistantly spaced points that span the domain. In addition, each table includes the reference where the competing method's solution error was found.

Although the TFC approach to solve differential equations has typically used ($N \in [100, 200$] and $m \in [20, 80]$) \cite{LDE,NDE,Selected,SVM}, in order to make a commensurable comparison with the techniques developed in Refs. \cite{R00,R03,R04,R05,R06,R07,R08,R09,R10,R11,R12,R13,R14,R15,R16,R17,R18}, a grid of $N = 11$ points was used and $m = 10$ basis terms were selected: one less than the number of points selected. The following sections introduced the six most commonly solved eighth-order differential equations in Refs. \cite{R00,R03,R04,R05,R06,R07,R08,R09,R10,R11,R12,R13,R14,R15,R16,R17,R18}, and compare the TFC solution accuracy with that of the most accurate solution from these references. Furthermore, to supplement the theory given in the previous sections, a step-by-step procedure for TFC is laid out for the first linear and non-linear problem.

\subsection{Linear Eighth-order Problems}
\subsubsection{Problem \#1}
Consider the linear eighth-order differential equation solved in Refs. \cite{R00,R12,R18}
\begin{equation*}
    y^{(8)}(x) - y(x) = -8e^{x} \quad x \in[0,1]
\end{equation*}
subject to
\begin{equation*}
\begin{aligned}
y(0) &= 1 \qquad \qquad y(1) &&= 0\\
y'(0) &= 0 \qquad \qquad y'(1) &&= -e\\
y''(0) &= -1 \qquad \quad y''(1) &&= -2e\\
y'''(0) &= -2 \qquad \quad y'''(1) &&= -3e\\
\end{aligned}
\end{equation*}
which has the exact solution $y(x) = (1-x)e^{x}$.

From Eqs. \eqref{eq:general_ce_a_b} and \eqref{eq:general_ce_a_b_derv}, we are shown that the estimated solution and eighth order derivative take on the forms,
\begin{equation*}
    y(x,\B{\xi}) = \B{a}(x)\T \B{\xi} + b(x)
\end{equation*}
and 
\begin{equation*}
    y(x,\B{\xi}) = \B{a}^{(8)}(x)\T \B{\xi} + b^{(8)}(x)
\end{equation*}
Thus, using TFC, the differential equation can be re-written as,
\begin{equation*}
    \tilde{F}(x,\B{\xi}) = \B{a}^{(8)}(x)\T \B{\xi} + b^{(8)}(x) - \B{a}(x)\T \B{\xi} - b(x) + 8e^{x} = 0.
\end{equation*}
Discretizing the problem into points, $x_k$ where $k \in[0,N]$, and collecting terms yields,
\begin{equation*}
    \mathbb{A}(\B{x})\B{\xi} + \B{b}(\B{x}) = 0,
\end{equation*}
where,
\begin{equation*}
    \mathbb{A} = \begin{Bmatrix} 
    \Big\{ \B{a}^{(8)}(x_0) - \B{a}(x_0) \Big\}\T \\ 
    \vdots \\ 
    \Big\{ \B{a}^{(8)}(x_k) - \B{a}(x_k) \Big\}\T  \\ 
    \vdots \\ 
    \Big\{ \B{a}^{(8)}(x_N) - \B{a}(x_N) \Big\}\T  \end{Bmatrix}
    \quad \text{and} \quad
    \B{b} = \begin{bmatrix} 
    b^{(8)}(x_0) - b(x_0) + 8e^{x_0} \\ \vdots \\
    b^{(8)}(x_k) - b(x_k) + 8e^{x_k} \\ \vdots \\ 
    b^{(8)}(x_N) - b(x_N) + 8e^{x_N} \end{bmatrix}.
\end{equation*}
This system can be solved by least-squares to yield the unknown coefficients, $\B{\xi}$, which can then be substituted back into the constrained expression to give the TFC estimate of the solution.

Table \ref{tab:p1} shows the absolute error of the TFC solution and the solution from Ref. \cite{R18} at each of the 11 points.
\begin{table}[H]
\begin{center}
\begin{tabular}{ccc} 
\toprule
{x} & {TFC Absolute Error} & {Ref. \cite{R18} Absolute Error}\\ \bottomrule \midrule
{0} & {0} & {0} \\
{0.1} & {2.2204e-16}  & {6.3e-11}\\
{0.2} & {1.1102e-16}  & {6.5e-10}\\
{0.3} & {1.1102e-16} & {2.0e-09}\\
{0.4} & {1.1102e-16}  & {3.3e-09}\\
{0.5} & {1.1102e-16}  & {3.9e-09}\\
{0.6} & {6.6613e-16}  & {3.4e-09}\\
{0.7} & {2.7756e-15}  & {2.0e-09}\\
{0.8} & {3.8858e-15}  & {6.9e-10}\\
{0.9} & {8.4932e-15}  & {7.6e-11}\\
{1} & {0} & {0}\\
\bottomrule
\end{tabular}
\end{center}
\caption{Problem \#1: Absolute solution error.}
\label{tab:p1}
\end{table}
Table \ref{tab:p1} shows that the TFC solution error is orders of magnitude lower than the solution from Ref. \cite{R18} at all of the points in the domain, except the boundaries. At the boundaries, each of the methods has zero error, because each of the methods satisfies the boundary conditions exactly. 

\subsubsection{Problem \#2}
Consider the linear eighth-order differential equation solved in Refs. \cite{R00,R17,R18,R13,R14,R12}
\begin{equation*}
    y(x)^{(8)} + xy(x) = -e^{x}(48 +15x + x^3) \quad x \in[0,1]
\end{equation*}
subject to
\begin{equation*}
\begin{aligned}
y(0) &= 0 \qquad \qquad y(1) &&= 0\\
y'(0) &= 1 \qquad \qquad y'(1) &&= -e\\
y''(0) &= 0 \qquad \qquad y''(1) &&= -4e\\
y'''(0) &= -3 \qquad \quad y'''(1) &&= -9e\\
\end{aligned}
\end{equation*}
which has the exact solution $y(x) = x(1-x)e^{x}$.

Table \ref{tab:p2} shows the absolute error of the TFC solution and the solution from Ref. \cite{R12} at each of the 11 points. Reference \cite{R12} did not report the solution at $x=0.9$, so that entry in the table is labeled ``not reported.''
\begin{table}[H]
\begin{center}
\begin{tabular}{ccc} 
\toprule
{x} & {TFC Absolute Error} & {Ref. \cite{R12} Absolute Error}\\ \bottomrule \midrule
{0} & {0} & {0}  \\
{0.1} & {0} & {1.63e-10}  \\
{0.2} & {8.3267e-17} & {1.63e-09}  \\
{0.3} & {0} & {4.90e-09}  \\
{0.4} & {1.1102e-16} & {8.46e-09}  \\
{0.5} & {5.5511e-17} & {1.01e-08}  \\
{0.6} & {3.8858e-16} & {8.68e-09}  \\
{0.7} & {3.3307e-16} & {5.15e-09}  \\
{0.8} & {3.3307e-16} & {1.76e-09}  \\
{0.9} & {8.0769e-15} & {Not reported}  \\
{1} & {0} & {0}  \\\midrule
\bottomrule
\end{tabular}
\end{center}
\caption{Problem \#2: Absolute solution error.}
\label{tab:p2}
\end{table}
Table \ref{tab:p2} shows that the TFC solution error is orders of magnitude lower than the solution from Ref. \cite{R12} at all of the points in the domain, except the boundaries. At the boundaries, each of the methods has zero error, because each of the methods satisfies the boundary conditions exactly. 

\subsubsection{Problem \#3}
Consider the linear eighth-order differential equation solved in Refs. \cite{R00,R18,R12,R05}
\begin{equation*}
    y(x)^{(8)} - y(x) = -8(2x\cos(x) + 7\sin(x)) \quad x \in[0,1]
\end{equation*}
subject to
\begin{equation*}
\begin{aligned}
y(0) &= 0 \qquad \qquad \quad y(1) &&= 0\\
y'(0) &= -1 \qquad \qquad y'(1) &&= -2\sin(1)\\
y''(0) &= 0 \qquad \qquad \quad y''(1) &&= 4\cos(1) + 2\sin(1)\\
y'''(0) &= 7 \qquad \qquad \quad y'''(1) &&= -6\sin(1) + 6\cos(1)\\
\end{aligned}
\end{equation*}
which has the exact solution $y(x) = (x^2-1)\sin(x)$.

Table \ref{tab:p3} shows the absolute error of the TFC solution and the solution from Ref. \cite{R18} at each of the 11 points.
\begin{table}[H]
\begin{center}
\begin{tabular}{ccc} 
\toprule
{x} & {TFC Absolute Error} & {Ref. \cite{R18} Absolute Error}\\ \bottomrule \midrule
{0} & {0} & {0}  \\
{0.1} & {2.7756e-17} & {6.6e-12}  \\
{0.2} & {2.7756e-17} & {6.9e-11}  \\
{0.3} & {0} & {2.1e-10}  \\
{0.4} & {5.5511e-17} & {3.5e-10}  \\
{0.5} & {0} & {4.1e-10}  \\
{0.6} & {7.2164e-16} & {3.5e-10}  \\
{0.7} & {1.3323e-15} & {2.1e-10}  \\
{0.8} & {1.1102e-15} & {7.2e-11}  \\
{0.9} & {3.4417e-15} & {8.0e-12}  \\
{1} & {0} & {0}  \\\midrule
\bottomrule
\end{tabular}
\end{center}
\caption{Problem \#3: Absolute solution error.}
\label{tab:p3}
\end{table}
Table \ref{tab:p3} shows that the TFC solution error is orders of magnitude lower than the solution from Ref. \cite{R18} at all of the points in the domain, except the boundaries. At the boundaries, each of the methods has zero error, because each of the methods satisfies the boundary conditions exactly. 

\subsection{Nonlinear eighth-order problems}
\subsubsection{Problem \#4}
Consider the nonlinear eighth-order differential equation solved in Refs. \cite{R00,R17}
\begin{equation*}
    y(x)^{(8)} + y^{(3)}(x)\sin(y(x)) = e^{x}(1+\sin(e^x)) \quad x \in[0,1]
\end{equation*}
subject to
\begin{equation*}
\begin{aligned}
y(0) &= 1 \qquad \qquad y(1) &&= e\\
y'(0) &= 1 \qquad \qquad y'(1) &&= e\\
y''(0) &= 1 \qquad \qquad y''(1) &&= e\\
y'''(0) &= 1 \qquad \qquad y'''(1) &&= e\\
\end{aligned}
\end{equation*}
which has the exact solution $y(x) = e^{x}$.

From Eqs. \eqref{eq:general_ce_a_b} and \eqref{eq:general_ce_a_b_derv}, we are shown that the estimated solution, third order derivative, and eighth order derivative take on the forms,
\begin{align*}
    y(x,\B{\xi}) &= \B{a}(x)\T \B{\xi} + b(x)\\
    y(x,\B{\xi}) &= \B{a}^{(3)}(x)\T \B{\xi} + b^{(3)}(x)\\
    y(x,\B{\xi}) &= \B{a}^{(8)}(x)\T \B{\xi} + b^{(8)}(x).
\end{align*}
Thus, using TFC, the differential equation can be re-written as,
\begin{equation*}
    \tilde{F}(x,\B{\xi}) = \B{a}^{(8)}(x)\T \B{\xi} + b^{(8)}(x) +  \Big( \B{a}^{(3)}(x)\T \B{\xi} + b^{(3)}(x) \Big) \sin\Big(\B{a}(x)\T \B{\xi} + b(x)\Big) - e^x\Big(1+\sin(e^x)\Big)= 0.
\end{equation*}
Discretizing the problem into points, $x_k$ where $k \in[0,N]$, leads to the loss function
\begin{equation*}
    \mathcal{L}(\B{\xi}_i) = \begin{bmatrix}  \tilde{F}(x_0,\B{\xi}) \\ \vdots \\
    \tilde{F}(x_k,\B{\xi}) \\ \vdots \\  \tilde{F}(x_N,\B{\xi}) \end{bmatrix}_{\B{\xi} = \B{\xi}_i},
\end{equation*}
for some values of $\B{\xi}=\B{\xi}_i$. The Jacobian of the loss function with respect to $\B{\xi}_i$ is,
\footnotesize
\begin{equation*}
\begin{aligned}
    &\mathcal{J}(\B{\xi}_i) = \begin{bmatrix} \dfrac{\partial \tilde{F}(x_0,\B{\xi})}{\partial \B{\xi}} \\ \vdots \\
    \dfrac{\partial \tilde{F}(x_k,\B{\xi})}{\partial \B{\xi}} \\ \vdots \\ 
    \dfrac{\partial \tilde{F}(x_N,\B{\xi})}{\partial \B{\xi}} \end{bmatrix}_{\B{\xi} = \B{\xi}_i} 
    \\ &= \begin{bmatrix} \Big\{ \B{a}^{(8)}(x_0) +  \Big( \B{a}^{(3)}(x_0) \Big) \sin\Big(\B{a}(x_0)\T \B{\xi} + b(x_0)\Big) + \Big(\B{a}^{(3)}(x_0)\T \B{\xi} + b^{(3)}(x_0)\Big) \cos\Big(\B{a}(x_0)\T \B{\xi} + b(x_0)\Big)\B{a}(x_0)  \Big\}\T \\ 
    \vdots \\
    \Big\{ \B{a}^{(8)}(x_k) +  \Big( \B{a}^{(3)}(x_k) \Big) \sin\Big(\B{a}(x_k)\T \B{\xi} + b(x_k)\Big) + \Big(\B{a}^{(3)}(x_k)\T \B{\xi} + b^{(3)}(x_k)\Big) \cos\Big(\B{a}(x_k)\T \B{\xi} + b(x_k)\Big)\B{a}(x_k) \Big\}\T\\ 
    \vdots \\
    \Big\{ \B{a}^{(8)}(x_N) +  \Big( \B{a}^{(3)}(x_N) \Big) \sin\Big(\B{a}(x_N)\T \B{\xi} + b(x_N)\Big) + \Big(\B{a}^{(3)}(x_N)\T \B{\xi} + b^{(3)}(x_N)\Big) \cos\Big(\B{a}(x_N)\T \B{\xi} + b(x_N)\Big)\B{a}(x_N) \Big\}\T
    \end{bmatrix}_{\B{\xi} = \B{\xi}_i} 
\end{aligned}
\end{equation*}
\normalsize

This system can be solved by an iterative least-squares method as shown in Section \ref{sec:SolveBVP} to yield the unknown coefficients, $\B{\xi}$, which can then be substituted back into the constrained expression to give the TFC estimate of the solution.

Table \ref{tab:p4} shows the absolute error of the TFC solution, which converged in 3 iterations, and the solution from Ref. \cite{R17} at each of the 11 points.
\begin{table}[H]
\begin{center}
\begin{tabular}{ccc} 
\toprule
{x} & {TFC Absolute Error} & {Ref. \cite{R17} Absolute Error}\\ \bottomrule \midrule
{0} & {0} & {0}  \\
{0.1} & {2.2204e-16} & {2.503395E-06}  \\
{0.2} & {0} & {8.940697E-06}  \\
{0.3} & {2.2204e-16} & {1.561642E-05}  \\
{0.4} & {4.4409e-16} & {1.823902E-05}  \\
{0.5} & {2.2204e-16} & {8.821487E-06}  \\
{0.6} & {6.6613e-16} & {7.510185E-06}  \\
{0.7} & {3.5527e-15} & {1.883507E-05}  \\
{0.8} & {7.5495e-15} & {1.931190E-05}  \\
{0.9} & {1.0214e-14} & {1.168251E-05}  \\
{1} & {0} & {0}  \\\midrule
\bottomrule
\end{tabular}
\end{center}
\caption{Problem \#4: Absolute solution error.}
\label{tab:p4}
\end{table}
Table \ref{tab:p4} shows that the TFC solution error is orders of magnitude lower than the solution from Ref. \cite{R17} at all of the points in the domain, except the boundaries. At the boundaries, each of the methods has zero error, because each of the methods satisfies the boundary conditions exactly. 

\subsubsection{Problem \#5}
Consider the nonlinear eighth-order differential equation solved in Refs. \cite{R00,R13,R14,R15,R17}
\begin{equation*}
    y(x)^{(8)} = 7! \Bigg(e^{-8y(x)} - \dfrac{2}{(1+x)^8}\Bigg) \quad x \in[0,e^{1/2} - 1]
\end{equation*}
subject to
\begin{equation*}
\begin{aligned}
y(0) &= 0 \qquad \qquad \qquad y(e^{1/2} - 1) &&= 1/2\\
y'(0) &= 1 \qquad \qquad \qquad y'(e^{1/2} - 1) &&= e^{-1/2}\\
y''(0) &= -1 \qquad \qquad \quad y''(e^{1/2} - 1) &&= -e^{-1}\\
y'''(0) &= 2 \qquad \qquad \qquad y'''(e^{1/2} - 1) &&= 2e^{-3/2}\\
\end{aligned}
\end{equation*}
which has the exact solution $y(x) = \ln(1+x)$.

Table \ref{tab:p5} shows the absolute error of the TFC solution, which converged in 2 iterations, and the solution from Ref. \cite{R14} at each of the 11 points.
\begin{table}[H]
\begin{center}
\begin{tabular}{ccc} 
\toprule
{x} & {TFC Absolute Error}  & Ref. {\cite{R14} Absolute Error}\\ \bottomrule \midrule
{0} & {0} & {0}  \\
{0.1} & {1.5266e-16} & {2.01e-07}  \\
{0.2} & {1.5821e-15} & {4.54e-07}  \\
{0.3} & {7.0083e-14} & {1.52e-06}  \\
{0.4} & {2.5846e-13} & {4.07e-06}  \\
{0.5} & {3.233e-13} & {6.71e-06} \\
{0.6} & {1.3139e-13} & {9.06e-06}\\
{0.7} & {2.1261e-14} & {1.00e-05}\\
{0.8} & {2.0539e-14} & {5.45-06}\\
{0.9} & {3.3307e-16} & {2.59e-06}\\
{1} & {0} & {0}  \\\midrule
\bottomrule
\end{tabular}
\end{center}
\caption{Problem \#5: Absolute solution error.}
\label{tab:p5}
\end{table}
Table \ref{tab:p5} shows that the TFC solution error is orders of magnitude lower than the solution from Ref. \cite{R14} at all of the points in the domain, except the boundaries. At the boundaries, each of the methods has zero error, because each of the methods satisfies the boundary conditions exactly. 
 
\subsubsection{Problem \#6}
Consider the nonlinear eighth-order differential equation solved in Refs. \cite{R00,R14,R17,R18}
\begin{equation*}
    y(x)^{(8)} + e^{-x} y^2(x) = e^{-x} + e^{-3x} \quad x \in[0,1]
\end{equation*}
subject to
\begin{equation*}
\begin{aligned}
y(0) &= 1 \qquad \qquad \qquad y(e^{1/2} - 1) &&= e^{-1}\\
y'(0) &= -1 \qquad \qquad \quad y'(e^{1/2} - 1) &&= -e^{-1}\\
y''(0) &= 1 \qquad \qquad \qquad y''(e^{1/2} - 1) &&= e^{-1}\\
y'''(0) &= -1 \qquad \qquad \quad y'''(e^{1/2} - 1) &&= -e^{-1}\\
\end{aligned}
\end{equation*}
which has the exact solution $y(x) = e^{-x}$.

Table \ref{tab:p6} shows the absolute error of the TFC solution, which converged in 2 iterations, and the solution from Ref. \cite{R18} at each of the 11 points.
\begin{table}[H]
\begin{center}
\begin{tabular}{ccc} 
\toprule
{x} & {TFC Absolute Error} & {Ref. \cite{R18} Absolute Error}\\ \bottomrule \midrule
{0} & {0} & {0}  \\
{0.1} & {1.1102e-16} & {2.9e-12}  \\
{0.2} & {1.1102e-16} & {2.7e-11}  \\
{0.3} & {0} & {7.6e-11}  \\
{0.4} & {0} & {1.3e-10}  \\
{0.5} & {1.1102e-16} & {1.5e-10}  \\
{0.6} & {1.1102e-16} & {1.3e-10}  \\
{0.7} & {2.2204e-16} & {7.6e-11}  \\
{0.8} & {3.2196e-15} & {2.5e-11}  \\
{0.9} & {9.992e-16} & {2.4e-12}  \\
{1} & {0} & {0}  \\\midrule
\bottomrule
\end{tabular}
\end{center}
\caption{Problem \#6: Absolute solution error.}
\label{tab:p6}
\end{table}
Table \ref{tab:p6} shows that the TFC solution error is orders of magnitude lower than the solution from Ref. \cite{R18} at all of the points in the domain, except the boundaries. At the boundaries, each of the methods has zero error, because each of the methods satisfies the boundary conditions exactly. 

\section{Accuracy of the Derivatives}\label{sec:derv}
The previous section compared the absolute error of TFC at a number of points along the domain with the absolute error of previous methods. This section discusses the accuracy of the derivatives when using TFC. One of the major advantages of TFC compared to other methods is the estimated solution is analytical. As a result, further manipulation of the estimated solution is easily achieved, for example, taking derivatives. Let's use problem \#5 as an example. 

Table \ref{tab:DerivativeError} shows the mean absolute error of $y$ and its derivatives up to order eight. The second column of Table \ref{tab:DerivativeError} used 10 basis functions to compute the solution and 11 equidistant points to compute the error, while the third column used 30 basis functions to compute the solution and 100 equidistant points to compute the error. 

\begin{table}[H]
\begin{center}
\begin{tabular}{ccc} 
\toprule
{Function} & \makecell{Mean Absolute Error:\\10 Basis Functions} & \makecell{Mean Absolute Error:\\30 Basis Functions} \\ \bottomrule \midrule
{$y$} &                    {7.5585e-14} & {9.6866e-16}\\
{$y^\prime$} &             {1.0534e-12} & {7.5884e-15}\\
{$y^{\prime\prime}$} &     {2.0202e-11} & {5.0360e-14} \\
{$y^{(3)}$} &              {4.9228e-10} & {4.0456e-13} \\
{$y^{(4)}$} &              {1.3318e-08} & {2.8079e-12} \\
{$y^{(5)}$} &              {3.8469e-07} & {1.3927e-11} \\
{$y^{(6)}$} &              {1.3150e-05} & {5.5250e-11} \\
{$y^{(7)}$} &              {3.9359e-04} & {2.0221e-10}\\
{$y^{(8)}$} &              {1.9399e-02} & {1.4765e-12}\\ \midrule
\bottomrule
\end{tabular}
\end{center}
\caption{Mean absolute error of all derivatives for Problem \#5}
\label{tab:DerivativeError}
\end{table}

If enough Chebyshev orthogonal polynomials are used in the free function to estimate the solution, the error in subsequent derivatives should increase by an order of magnitude or less.  Table \ref{tab:DerivativeError} shows that when 10 basis functions are used, the error increases as the order of the derivative increases. In this case, there were not enough Chebyshev orthogonal polynomials used, as indicated by the large mean error in the eighth derivative. In other words, the number of basis functions used was not nearly enough to accurately estimate the solution of the eighth derivative. 

When 30 basis functions are used, the mean error increases as the order of the derivative increases, until the eighth derivative is reached. In derivatives one through seven, the mean error increases approximately an order of magnitude or less when compared to the previous derivative. However, the eighth derivative has less error than the seventh derivative, because the eighth derivative shows up in the differential equation, and thus in the residual. Hence, the eighth derivative is directly affected when computing the solution, whereas the other derivatives are not.  

In problem \#5, the differential equation only contains the function and the eighth derivative. As a different example, consider the problem solved in \cite{R16},
\begin{equation*}
\begin{aligned}
    y(x)^{(8)} + y(x)^{(7)} + 2y(x)^{(6)} + 2y(x)^{(5)} &+ 2y(x)^{(4)} + 2y(x)^{(3)} + 2y(x)^{(2)} + y(x)^{(1)} + y(x) \\ &= 14\cos(x) - 16\sin(x) - 4x\sin(x) \quad x \in[0,1]
\end{aligned}
\end{equation*}
subject to
\begin{equation*}
\begin{aligned}
y(0) &= 0 \qquad \qquad \quad y(1) &&= 0\\
y'(0) &= -1 \qquad \quad \quad y'(1) &&= 2 \sin(1)\\
y''(0) &= 0 \qquad \qquad \quad y''(1) &&= 4\cos(1) + 2\sin(1)\\
y'''(0) &= 7 \qquad \qquad \quad y'''(1) &&= 6\cos(1) - 6\sin(1)\\
\end{aligned}
\end{equation*}
which has the exact solution $y(x) = (x^2 - 1)\sin(x)$. From hereon, we shall refer to this problem as problem \#7.

Table \ref{tab:ExtraProblem} shows the mean absolute error of $y$ and its derivatives up to order eight for problem \#7. The second column of Table \ref{tab:ExtraProblem} used 10 basis functions to compute the solution and 11 points to compute the error, while the third column used 30 basis functions to compute the solution and 100 points to compute the error. 

\begin{table}[H]
\begin{center}
\begin{tabular}{ccc} 
\toprule
{Function} & \makecell{Mean Absolute Error:\\10 Basis Functions} & \makecell{Mean Absolute Error:\\30 Basis Functions} \\ \bottomrule \midrule
{$y$} &                    {4.8255e-15} &{5.8919e-15}\\
{$y^\prime$} &             {8.3368e-15} &{9.9755e-15}\\
{$y^{\prime\prime}$} &     {7.1054e-14} &{5.9525e-14} \\
{$y^{(3)}$} &              {4.8760e-13} &{4.8352e-13} \\
{$y^{(4)}$} &              {1.5118e-12} &{1.6443e-12} \\
{$y^{(5)}$} &              {4.1244e-12} &{4.4041e-12} \\
{$y^{(6)}$} &              {3.1934e-12} &{3.2911e-12} \\
{$y^{(7)}$} &              {8.0532e-12} &{7.3956e-12}\\
{$y^{(8)}$} &              {8.1927e-11} &{8.7722e-12}\\ \midrule
\bottomrule
\end{tabular}
\end{center}
\caption{Mean absolute error of all derivatives for Problem \#7}
\label{tab:ExtraProblem}
\end{table}

Table \ref{tab:ExtraProblem} shows that when all derivatives are included in the differential equation, the anomalous decrease in mean error as subsequent derivatives are taken disappears (i.e. the mean solution error from derivative seven to derivative eight increases as expected). 

\section{Conclusions}

This paper explores the application of the techniques developed in \cite{U-TFC,LDE,NDE} to the solution of high-order differential equations, namely eighth-order BVPs. In all the problems presented, which span the publications \cite{R00,R03,R04,R05,R06,R07,R08,R09,R10,R11,R12,R13,R14,R15,R16,R17,R18}, the solution accuracy ranges from $\mathcal{O}(10^{-13} - 10^{-16})$. These results are similar to the results obtained in earlier studies on first- and second-order linear \cite{LDE} and nonlinear \cite{NDE} differential equations. This application to higher-order systems further highlights the power and robustness of this technique.

In Section \ref{sec:init}, a discussion of the initialization of the TFC approach for nonlinear differential equations was provided. In the initialization, the coefficient vector $\B{\xi}$ is set to zero, which implies that $g(x) = 0$. It was found that this still solved the differential equation with an accuracy on the order of $\mathcal{O}(10^{-7})$. This can be explained by an equation first presented in the seminal paper on TFC \cite{U-TFC}. In this paper an equation (Eq. (9) in Ref. \cite{U-TFC}) is presented which describes the general expression for the interpolating expression for the function and its first $n$ derivatives. This equation simplifies to the expression of a Taylor series when $g(x) = 0$. While this is not the exact case for the boundary-value constraints in this paper, the nature of the aforementioned equation points to the fact that when $g(x) = 0$, the \ce\ derived in this paper acts as a Taylor series approximation for two points. The relationship between the TFC method and Taylor series has yet to be explored, and the extension of this Taylor series-like expansion about $n$ points will be the focus of future work. 

Section \ref{sec:derv} discussed the accuracy of the derivatives of the estimated solution. The solution accuracy was reduced with each subsequent derivative, but overall the accuracy of the final derivatives only lost a few orders of magnitude, resulting in an overall error on the order of $\mathcal{O}(10^{-10} - 10^{-12})$, provided that enough basis functions were used when estimating the solution. In addition, it was shown that the accuracy of a given derivative depends, in part, on whether it explicitly shows up in the differential equation. 
\vspace{6pt} 

 
\appendix
\section{Support Functions for General Points \texorpdfstring{$x_i$}{xi} and \texorpdfstring{$x_f$}{xf}}\label{app:SupportFunctions}
This appendix shows the switching functions for the constrained expression shown in Eq. \eqref{eq:ce_solved} for a general domain $x \in [x_i, x_f]$. The switching functions are:

\small
\begin{align*}
\beta_1(x) &= \frac{\left(x-x_f\right){}^4 \left(-7 \left(4 x x_f+x_f^2+10 x^2\right) x_i+21 \left(x_f+4 x\right) x_i^2+10 x^2 x_f+4 x x_f^2+x_f^3-35 x_i^3+20 x^3\right)}{\left(x_f-x_i\right){}^7}\\
\beta_2(x) &= -\frac{\left(x-x_i\right){}^4 \left(-7 x_f \left(4 x x_i+x_i^2+10 x^2\right)+21 x_f^2 \left(x_i+4 x\right)-35 x_f^3+10 x^2 x_i+4 x x_i^2+x_i^3+20 x^3\right)}{\left(x_f-x_i\right){}^7}\\
\beta_3(x) &= \frac{\left(x-x_f\right){}^4 \left(x-x_i\right) \left(-6 \left(x_f+4 x\right) x_i+4 x x_f+x_f^2+15 x_i^2+10 x^2\right)}{\left(x_f-x_i\right){}^6}\\
\beta_4(x) &= \frac{\left(x-x_f\right) \left(x-x_i\right){}^4 \left(-6 x_f \left(x_i+4 x\right)+15 x_f^2+4 x x_i+x_i^2+10 x^2\right)}{\left(x_f-x_i\right){}^6}\\
\beta_5(x) &= \frac{\left(x-x_f\right){}^4 \left(x-x_i\right){}^2 \left(x_f-5 x_i+4 x\right)}{2 \left(x_f-x_i\right){}^5}\\
\beta_6(x) &= -\frac{\left(x-x_f\right){}^2 \left(x-x_i\right){}^4 \left(-5 x_f+x_i+4 x\right)}{2 \left(x_f-x_i\right){}^5}\\
\beta_7(x) &= \frac{\left(x-x_f\right){}^4 \left(x-x_i\right){}^3}{6 \left(x_f-x_i\right){}^4}\\
\beta_8(x) &= \frac{\left(x-x_f\right){}^3 \left(x-x_i\right){}^4}{6 \left(x_f-x_i\right){}^4}
\end{align*}
\normalsize

\bibliographystyle{unsrt} 
\bibliography{mybib} 
\end{document}